\documentclass[11pt]{amsart2000}
\usepackage{epsfig}
\newtheorem{lem}{Lemma}[section]
\newtheorem{prop}[lem]{Proposition}
\newtheorem{thm}[lem]{Theorem}
\newtheorem{cor}[lem]{Corollary}
\theoremstyle{definition}
\newtheorem{defn}[lem]{Definition}
\newtheorem{rem}[lem]{Remark}
\newtheorem{exam}[lem]{Example}

\begin{document}
\setlength{\unitlength}{0.01in}
\linethickness{0.01in}
\begin{center}
\begin{picture}(474,66)(0,0)
\multiput(0,66)(1,0){40}{\line(0,-1){24}}
\multiput(43,65)(1,-1){24}{\line(0,-1){40}}
\multiput(1,39)(1,-1){40}{\line(1,0){24}}
\multiput(70,2)(1,1){24}{\line(0,1){40}}
\multiput(72,0)(1,1){24}{\line(1,0){40}}
\multiput(97,66)(1,0){40}{\line(0,-1){40}}
\put(143,66){\makebox(0,0)[tl]{\footnotesize Proceedings of the Ninth Prague Topological Symposium}}
\put(143,50){\makebox(0,0)[tl]{\footnotesize Contributed papers from the symposium held in}}
\put(143,34){\makebox(0,0)[tl]{\footnotesize Prague, Czech Republic, August 19--25, 2001}}
\end{picture}
\end{center}
\vspace{0.25in}
\setcounter{page}{75}
\title{On finite $T_0$ topological spaces}
\author{A. El-Fattah El-Atik}
\address{Ernst-Moritz-Arndt-Universit\"at,
Institut f\"ur Mathematik und Informatik, Jahnstr. 15a, 17487
Greifswald, Germany} 
\email{elatik@mail.uni-greifswald.de}
\address{From November 20,2001 the new address will be\\
Department of Mathematics, Faculty of Science,
Tanta University, Tanta, Egypt}
\email{aatik@dec1.tanta.eun.eg}
\author{M. E. Abd El-Monsef}
\address{Department of Mathematics, Faculty of Science,
Tanta University, Tanta, Egypt}
\author{E. I. Lashin}
\address{Department of Mathematics, Faculty of Engineering,
Tanta University, Tanta, Egypt}
\subjclass[2000]{Primary: 54B10, 54D30; Secondary: 54A05, 54G99}
\keywords{$T_0$ spaces, finite $T_0$ spaces, minimal neighborhoods, 
partially ordered sets, multifunctions, topological dimension, weak 
continuity}
\begin{abstract}
Finite topological spaces became much more essential in topology, with the
development of computer science. 
The task of this paper is to study and investigate some properties of such
spaces with the existence of an ordered relation between their minimal
neighborhoods. 
We introduce notations and elementary facts known as Alexandroff
space~\cite{arenas1, arenas2, adg}. 
The family of minimal neighborhoods forms a unique minimal base. 
We consider $T_0$ spaces. 
We give a link between finite $T_0$ spaces and the related partial order. 
Finally, we study some properties of multifunctions and their
relationships with connected ordered topological spaces.
\end{abstract}
\thanks{A. El-Fattah El-Atik, M. E. Abd El-Monsef and E. I. Lashin,
{\em On finite $T_0$ topological spaces},
Proceedings of the Ninth Prague Topological Symposium, (Prague, 2001),
pp.~75--90, Topology Atlas, Toronto, 2002}
\maketitle

\section{Introduction and Preliminaries}

Finite spaces were first studied by P.A.~Alexandroff in 1937 
in~\cite{alex}. 
Actually, finite spaces had been more earlier investigated by many authors
under the name of simplicial complexes. 
There were several other contributions by Flachsmeyer in 1961 \cite{flac},
Stong in 1966 \cite{ston} and L. Lotz in 1970 \cite{lot}. 
Rinow \cite{rin} in his book discussed some properties of finite spaces. 
However, the subject has never been considered as a main field of
topology.

With the progress of computer technology, finite spaces have become more
important. 
Herman in 1990 \cite{her}, Khalimsky and et.~al.\ in 1990 \cite{khal},
Kong and Kopperman in 1991 \cite{kon} have been applied them to model the
computer screen. In this paper we focus on finite spaces with
order. 
The main importance of our study is to offer a new formulations for some
topological operators in general topology such as interior, closure,
boundary and exterior operators. 
We present and study comparisons between some topological properties in
the case of finite spaces. 
In what follows, by $X$ we mean always mean a finite $T_0$ space. 
For each $A\subset X$, the closure (resp.\ interior, exterior,
boundary) of $A$ will be denoted by $\overline{A}$ 
(resp.\ $\mbox{ int }(A)$, $\mbox{ ext }(A)$, $\partial A$).

For each point $x$ in a space $X$, there is a smallest neighborhood which
is contained in each other neighborhood of $x$. 
For each $x\in X$, let
\[U_x = \bigcap \{V: V \mbox{ is an open set containing }x \}\] Clearly
$U_x$ is the smallest open set containing $x$ since $X$ is finite.\\

Alexandroff spaces are the topological spaces in which each element is
contained in a smallest open set or equivalently the spaces where
arbitrary intersections of open sets are open. 
It is clear that all finite
spaces are locally finite and all locally finite spaces are Alexandroff.

\begin{lem} \label{l:baseft}
The class $\mathcal{U} = \{U_x: x\in X\}$ is a base for a finite space $(X, \tau)$. 
Each base of $\tau$ contains $\mathcal{U} $.
\end{lem}

Notice that if $X$ is Alexandroff, then $X$ is $T_1$ if and only if 
$U_x= \{x\}$. 
It follows that $X$ is discrete and so every point is an isolated point.

\begin{rem}\label{r:twominnbd}
Observe that if $x$ and $y$ are two points in a space $X$, then $y\in U_x$
if and only if $U_y \subseteq U_x$.
\end{rem}

\begin{defn}\label{d:relation}\cite{alek}
For two points $x, y\in X$, ~$y\geq x$ if ~$U_y\subseteq U_x$.
\end{defn}

\begin{rem}
From Definition~\ref{d:relation}, the relation $\geq $ is reflexive and
transitive since $\subseteq$ is so.
\end{rem}

\begin{prop}\label{p:rewclosure}
In a space $X$, $y\geq x$ if and only if $x\in \overline{\{y\}}$. 
In such case, $x$ is said to be identified with $y$.
\end{prop}

\begin{proof}
Let $y\geq x$ and $y\neq x$. 
Then $y\in U_x$ which is the smallest open set containing $x$. 
Then for any open set $G$ containing $x$, we have 
$(G\setminus \{x\})\cap \{y\}\neq \phi$. 
This means $x$ is an accumulation point of $y$. 
Therefore $x\in \overline{\{y\}}$. 
Conversely, let $x\in \overline{\{y\}}$. 
Then $G\cap \{y\}\neq \phi$ for every open sets $G$ containing $x$. 
So $y\in G$ for every open set $G$. 
Take $G= U_x$. 
By Remark~\ref{r:twominnbd}, we get $U_y\subseteq U_x$. 
This shows that $y\geq x$.
\end{proof}

The following example will be used throughout the paper. 

\begin{exam}[The class of topologies on a set with three points]
\label{ex:thptnbd}
We give a list of all topologies on the set $X= \{x, y, z\}$ up to
homeomorphic topologies. 
There are 9 topologies, with the following $U_x$, $U_y$ and $U_z$:
\begin{itemize}
\item[\textbf{\large $\tau_1$}] 
All three points are isolated. This is the discrete topology and $U_x=
\{x\}$, $U_y=\{y\}$, $U_z= \{z\}$.
\end{itemize}
Now, we assume that the topology has only two isolated points ,
say, $U_x=\{x\}$ and $U_y= \{y\}$. We have two cases concerning
$U_z$:
\begin{itemize}
\item[\textbf{\large $\tau_2 $}]
The neighborhood of $U_z$ has two points. We can assume $U_z= \{x, z\}$.
\item[\textbf{\large $\tau_3 $}] 
The neighborhood of $U_z$ has three points, $U_z= X$.
\end{itemize} 
Next, we consider the case that the topology has only one isolated
point, say, $U_x= \{x\}$. Again, we have to distinguish between
different cases concerning $U_y$ and $U_z$:
\begin{itemize}
\item[\textbf{\large $\tau_4$}] 
Both $U_y$ and $U_z$ have three points. Then $U_y= U_z= X$.
\item[\textbf{\large $\tau_5$}] 
One neighborhood, say $U_z$, has three points and $U_y$ has two
points. Then $U_y= \{y, z\}$ is not possible, so $U_y= \{x, y\}$.
\item[\textbf{\large $\tau_6$}] 
Both neighborhoods have two points and are equal $U_y= U_z= \{y, z\}$.
\item[\textbf{\large $\tau_7 $}] 
$U_x$ and $U_y$ have two points and are different. Then they can not be
$\{y, z\}$. Thus $U_y= \{x, y\}$ and $U_z= \{x, z\}$.
\end{itemize}
Finally, we have two cases without isolated points:
\begin{itemize}
\item[\textbf{\large $\tau_8$}] 
There is a neighborhood with two points, say, $U_x= \{x, y\}$. Then the
neighborhood of $z$ must have three points. Thus $U_x= U_y= \{x, y\}$ and
$U_z= X$.
\item[\textbf{\large $\tau_9$}]
All neighborhoods have three points, then $U_x= U_y= U_z= X$.
\end{itemize}
Figure~\ref{f:figps} shows this class of topologies with its minimal
neighborhoods. 
The smallest ellipse refers to a singleton, the middle is a two points
neighborhood and the biggest one is the whole space $X$.

\begin{figure}
\centering \epsfxsize100mm \epsfysize120mm
\centerline{\epsffile{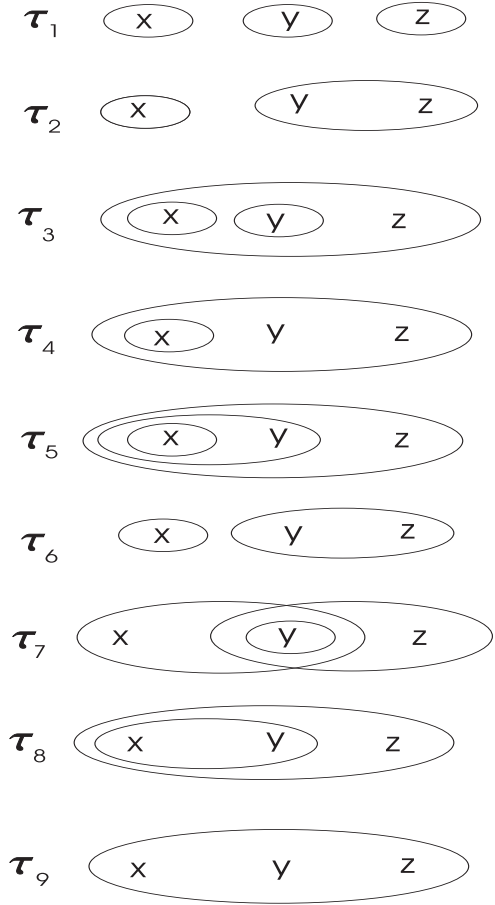}} 
\caption{Topological spaces with minimal neighborhoods.}\label{f:figps}
\end{figure}
\end{exam}

\section{$T_0$ Properties and Associated Partial Order}

In this section, we investigate some properties of $T_0$ spaces
such that for two distinct points $x$, $y$ in $X$, $x\geq y$ and
$y\geq x$ only true if $U_x= U_y$. If $\mathcal{U} (x)$ denote the
neighborhood system of $x$. Recall that a space is $T_0$ if and
only if $\mathcal{U} (x)\neq \mathcal{U} (y)$ for $x\neq y$.

\begin{rem}
The topologies $\tau_i$ for $i\in \{1, 2, 3, 5, 7\}$, in
Example~\ref{ex:thptnbd}, are $T_0$. Any two elements in
Figure~\ref{f:figps} with the same neighborhood is not $T_0$ such
as $\tau_j$ for $j\in \{4, 6, 8, 9\}$.
\end{rem}

\begin{defn}\label{d:twonbdsys}
For points $x, y$ in a topological space $X$, $y\geq x$ if $\mathcal{U} 
(y)\supseteq \mathcal{U} (x)$. Notice that $\supseteq$ means more numbers
of neighborhoods.
\end{defn}

\begin{rem}
The relation $\geq$ is reflexive and transitive, since this is
true for $\supseteq$. In a $T_0$ space, it is also antisymmetric.
Then $\geq$ is a partial order in $T_0$ spaces.
\end{rem}

\begin{defn}\label{d:openinparor}
Let $(X, \leq)$ be a partially ordered set and $U\subseteq X$. We
say that $U$ is open if whenever $x\in U$ and $y\geq x$ it is also
the case that $y\in U$.
\end{defn}

In the following result we denote by the set of open subsets of
$X$ by $\mathcal{U} _ \leq $.

\begin{prop}
If $(X, \leq)$ is a partially ordered set, then $(X, \mathcal{U} _{\leq})$
is a finite $T_0$ space.
\end{prop}
\begin{proof}
Clearly, $X$ and $\phi$ are elements in $\mathcal{U} _\leq$. 
Let $U_i \in \mathcal{U} _\leq \mbox{ for every } i\in I$. 
For any $x\in \bigcup\limits_{i\in I}U_i$ and $y\geq x$, there is 
$i_0\in I$ such that $x\in U_{i_0}$. 
By openness of $U_{i_0}$, we have $y\in U_{i_0}$ and 
$y\in \bigcup\limits_{i\in I}U_i$. 
Therefore
$\bigcup\limits_{i\in I}U_i\in \mathcal{U} _\leq$. 
Also, if $A$ and $B$ are elements in $\mathcal{U} _\leq$, then 
$A\cap B \in \mathcal{U} _\leq$. 
To show $T_0$, consider two distinct elements $x, y \in X$. 
Clearly, $x\in U_x$ and $y\in U_y$. 
If $y\not\in U_x$, the proof is complete. 
If $y\in U_x$, $y\geq x$, by the antisymmetry of $\leq$, $x\not\geq y$
and so $x\notin U_y$ which also completes the proof.
\end{proof}

\begin{rem}
The closure of any singleton $p$ of a finite $T_0$ space $(X,
\mathcal{U} _\leq)$ has the form $\overline {\{p\}}= \{x : x\leq p\}$.
\end{rem}

\begin{rem}
The order relation $\leq$ is not useful in $T_1$ spaces, because
there are no comparable elements. In other words, since for every
two distinct points $x$ and $y$ in a $T_1$ space,
$\overline{\{x\}}\neq \overline{\{y\}}$, then the relation $x\leq
y$ is never satisfied. Thus $\leq$ is a good tool only for $T_0$
spaces which are not $T_1$.
\end{rem}

\begin{prop}\label{p:tswisolated}
A finite $T_0$ space contains an isolated point.
\end{prop}
\begin{proof}
Obvious.
\end{proof}

\begin{rem}
Finite spaces with isolated points need not be $T_0$. Each of
$\tau_4$ and $\tau_6$ in figure~\ref{f:figps} of
example~\ref{ex:thptnbd} has an isolated points, but not $T_0$.
\end{rem}

\begin{prop}\label{p:openisolated}
Every open set in a finite $T_0$ space contains an isolated point.
\end{prop}

\begin{proof}
Let $G$ be an arbitrary open set in $X$. Then $G$ is a finite open
$T_0$ subspace. By Proposition~\ref{p:tswisolated}, there is an
isolated point $x$ in $G$. Since $G$ is open in $X$, then $x$ must
be isolated point in $X$.
\end{proof}

\section{Some Topological Properties}

In this paper, some of topological operators which are well known
for topologists have new forms with respect to the Alexandroff's
notion of order.

\begin{prop} \label{p:upperlower}
In any finite space, the open (closed) points are the maximal
(minimal) elements. Therefore $U$ is an open set if it contains
its upper bounds. By the complement, $F$ is closed if it contains
its lower bounds.
\end{prop}
\begin{proof}
By Definition~\ref{d:relation} and Proposition~\ref{p:rewclosure},
the partial order $x\leq y$ for points $x, y$ in a finite $T_0$
space $X$ was defined by $y\in U_x$ or $x\in \overline {\{y\}}$.
This means that for open points $\{x\}= U_x$, there is no strictly
larger $y$. Thus the open points are the maximal elements in $X$.
Similarly, for a closed point $y$, if $U_x= \{x\}$, then $x\leq y$
implies $y= x$. Hence the closed points are the minimal elements.
\end{proof}

\begin{prop}\label{p:interiorset}
In a finite topological space $X$, the interior of a subset
$A\subset X$, $\mbox{ int }A$, is the set of all points $a\in A$
such that for every $b\geq a$ it is also true that $b\in A$. In
other words, the elements of $\mbox{ int }A$ are all upper bounds
of $X$ belonging to $A$.
\end{prop}
\begin{proof}
We prove that $\mbox{ int }A= \{a\in A: b\geq a \mbox{ implies }
b\in A\}$. For $``\subseteq''$, let $a\in \mbox{int}(A)$ that is
$U_a\subseteq A$ for $\mbox{ int }A$ is the greatest open set
contained in $A$. But $U_a= \{b: b\geq a\}$. For $``\supseteq''$, if
$b\in A$ for all $b\geq a$. That means $U_a \subseteq A$, then
$a\in \mbox{ int }A$.
\end{proof}

\begin{prop} \label{p:closureset}
For a subset $A$ in a finite topological space $X$, the closure of
$A$, $\overline A$, is the set of points $b\in X$ such that $b\leq
a$ for some $a\in A$. In other words, $\overline A$ is the set of
all lower bounds of points of $A$.
\end{prop}
\begin{proof}
In a finite topological space $X$, $\overline A=
\bigcup\limits_{a\in A}\overline {\{a\}}$. By the definition of
$\leq$, we have $\overline {\{a\}}= \{b: b\leq a\}$ and the claim.
\end{proof}

\begin{prop}\label{p:clint}
For a finite topological space $X$, the closure of an interior of
a subset $A$ of $X$, $\overline {\mbox{ int }A}$, is the set of
all $x$ in $X$ such that there exists a maximal element $a\in A$
with $a\geq x$.
\end{prop}
\begin{proof}
By Proposition~\ref{p:closureset}, $\overline{\mbox{ int }A}=
\bigcup\{b: b\leq a \mbox{ for some }a\in \mbox{ int }A \}$. Since
$a\in\mbox{int }A $, then all $c\geq a$ are also in $A$. Therefore
$\overline {\mbox{ int }A}= \bigcup\{b: b\leq a \mbox{ for }a\in A
\}$. Since we have a finite ordered set, there is a maximal $c\geq
a$. Hence $\overline {\mbox{ int }A}$ is the lower bounds of
maximal elements in $A$.
\end{proof}

\begin{prop}\label{p:intcl}
In a finite topological space $X$, the interior of a closure of
$A\subset X$,$\mbox{ int }\overline A$, is the points $x\in X$
such that all maximal elements $y\geq x$ belong to $A$.
\end{prop}
\begin{proof}
From~\ref{p:interiorset} and \ref{p:closureset}, we have
\begin{eqnarray}
\mbox{ int }\overline A &=& \{a\in \overline A: b\geq a \mbox{
implies }b\in \overline A \}\nonumber\\
&=& \{a\in \overline A: b\geq a \mbox{ implies that there is }c\in
A, c\geq b\}\nonumber\\
&=& \{a\in \overline A: \mbox{ all maximal elements }b\geq a
\mbox{ must belong to }A\}\nonumber
\end{eqnarray}
To prove the equality, if $b$ is maximal, then $b\in A$ since
there is no other $c$ in $A$. Hence $b\in \mbox{int}\overline A$.
On the other hand, if all elements in $A$ are maximal and $b\geq
a$, then there is a maximal $c\geq b$ with $c\in A$.
\end{proof}

In finite topological spaces, From Propositions~\ref{p:clint}~
and~\ref{p:intcl}, we have $\mbox{int}\overline A\subseteq
\overline {\mbox{ int }A}$. The following example shows that the
inverse inclusion is not always true.

\begin{exam}\label{ex:closint}
Consider the topology $\tau_3$ in Example~\ref{ex:thptnbd}. If $A=
\{y, z\}$, we get $\mbox{ int }{\overline A}= \{y\}$ and
$\overline{\mbox{ int }A}= \{y, z\}$. Therefore
$\mbox{int}\overline A \not\supseteq \overline {\mbox{ int } A}$.
\end{exam}

In 1961, Levine \cite{levine1} introduced the notion of semi-open
set in any topological space. A subset $A$ in a topological space
$X$ is called semi-open if and only if $A\subset \overline {\mbox{
int }A}$. In 1982, Mashhour et.~al.\ \cite{mash} defined the
concept of preopen set. A subset $A$ is preopen in $X$ if and only
if $A\subset \mbox{ int }{\overline A}$. There are no implications
between these two concepts. This means that semiopen sets need not
be preopen and conversely. In 1997, Abd El-Monsef et.~al.\
\cite{abdel} defined the concept of $\gamma$-open set as a union
of semi-open and preopen sets which is equivalent to $A\subset
\mbox{ int }{\overline A} \cup \overline {\mbox{ int }A}$, for any
subset $A$ of $X$. From the relation $\mbox{int}\overline
A\subseteq \overline {\mbox{ int }A}$, we obtain a new implication
between types of near openness in finite topological spaces.

\begin{rem}\label{r:semipre}
In a finite topological space $X$, the following are true:
\begin{enumerate}
\item[(i)] Every preopen set in $X$ is semi-open. The converse
may not be true as it is shown in Example~\ref{ex:closint}
\item[(ii)] Semiopen sets coincide with $\gamma$-open sets. 
The following diagram shows the relation between these notions in a
finite case.

$$\begin{array}{crccr}
\textbf{\large openness}&\Longrightarrow& \textbf {\large
preopenness}& \\
&&\Downarrow&& \\ 
\textbf{\large $\gamma $-openness }&\Longleftrightarrow
&\textbf{\large semi-openness}&&
\end{array}$$
\end{enumerate}
\end{rem}

The axiom $T_0$ is necessary for satisfying the previous
implications. The following examples show this fact.

\begin{exam}
Any subset of an indiscrete space is preopen but not semi-open,
since the indiscrete spaces are not $T_0$.
\end{exam}

\begin{prop}
In a finite space $X$, the boundary of a subset $A$ of $X$ is the
set $b\in X$ such that $b\leq a$ for some $a\in A$ and $b\leq c$
for some $c\notin A$.
\end{prop}
\begin{proof}
It follows from the definition, $\partial A= \overline A\cap
\overline {X\setminus A}$, and Proposition~\ref{p:closureset}.
\end{proof}

\begin{prop}
In a finite space $X$, the exterior of $A\subset X$ is the set of
points $x\in X$ such that for $b\geq x$ implies $b\notin A$.
\end{prop}
\begin{proof}
Since for any topological space $\mbox{ ext (A)}= X\setminus
{\overline A}= \mbox { int }(X\setminus A)$. By
Proposition~\ref{p:interiorset}, we have
\begin{eqnarray}
\mbox{ ext (A)}&=& \{ x\in X: b\geq x \mbox { implies }b\in
X\setminus A\}\nonumber\\
&=& \{x\in X: b\geq x \mbox { implies }b\notin A\}\nonumber
\end{eqnarray}
\end{proof}

\section{The Dimension for Finite Spaces}

The notion of a topological dimension has a sense for finite
spaces. Although we have only a finite number of points, all
finite dimensions are possible. We use the notion of inductive
dimension as stated by Engelking (\cite{eng}, Chapter 7). A
discrete space $X$ has dimension zero denoted by $\dim X= 1$,
since the neighborhoods have no boundaries. On the real line $R$,
the open intervals form a basis. The boundary of an interval
consist of two points which is a discrete space. Analogously, the
plane has dimension $\leq 2$ since a circumference has dimension
$\leq 1$. This is the typical situation in three-dimensional
spaces. This leads to the following definition.

\begin{defn}(\cite{cain}, p.195) \label{d:dim}
Let $X$ be a topological space. $\dim X= -1$ if and only if $X=
\phi$. Let $n$ be a positive integer and $\dim X\leq k$ be defined
for each $k\leq n-1$. Then $\dim X\leq n$ if $X$ has a base
$\beta$ such that $\dim\partial B\leq n-1$ for all $ B\in\beta $.
\end{defn}
A base $\beta$ of any finite space will be replaced by the set of
all minimal neighborhoods for each of its points.

\begin{prop}
For a finite space $X$, $\dim X\leq n$ if and only if the minimal
base $\mathcal{U} $ fulfils $\dim \partial U_x\leq n-1$ for each $U_x\in
\mathcal{U} $.
\end{prop}
\begin{proof}
If $\dim X\leq n$, then by definition~\ref{d:dim} it has a base
$\beta$ such that $\dim\partial B\leq n-1$ for all $ B\in\beta $.
By Lemma~\ref{l:baseft}, each $\beta$ contains $\mathcal{U} $. Then $\dim
\partial U_x\leq n-1$. The converse is obvious.
\end{proof}

\begin{prop}
A topological space $X$ has dimension zero if and only if $X$ has
a base of clopen sets.
\end{prop}
\begin{proof}
It is clear that $\partial B= \phi$ if and only if $B$ is clopen.
Dimension zero means existence of a base with $\partial B= \phi$.
This condition is equivalent to $B$ being clopen. Since $\partial
B$ is containing those $x$ which are accumulation points of $B$
and its complement. For clopen sets $B$, we have $\partial B=
\phi$. On the other hand, if $\partial B= \partial (X\setminus B)=
\phi$, then $B$ and $X\setminus B$ are closed and so $B$ is
clopen.
\end{proof}

\begin{prop}
Every finite $T_0$ space with base $\beta$ of clopen sets is
discrete.
\end{prop}
\begin{proof}
Since every zero dimensional $T_0$ space is $T_1$. then it is
closed as a consequence of the finiteness of the space.
\end{proof}

Finiteness and $T_0$ axiom are necessary. It can be satisfied from
the following examples.

\begin{exam}
The Cantor set $\{0, 1\}^{\infty}$ with product topology has a
base of clopen sets, but this space is not discrete.
\end{exam}

\begin{exam}
In Example~\ref{ex:thptnbd}, the topology $\tau_6$ has a base
$\{\{x\}, \{y, z\}\}$ with clopen sets. $X$ is not discrete,
because $X$ is not $T_0$.
\end{exam}

\begin{lem}\label{l:clopendiscrete}
Let $X= C\cup V$ such that every $c\in C$ is closed and $v\in V$
is open. Then each of $C$ and $V$ is a discrete subspace of $X$.
\end{lem}
\begin{proof}
Since $C$ is a finite closed subspace of $X$, then each $c\in C$
is a closed point in $C$. Then $C\setminus \{c\}$ is a finite
closed subset in $C$. So $\{c\}$ is an open point in $C$. This
means that $C$ is a discrete set. It is clear that $V$ is
discrete.
\end{proof}

The following Proposition describes the finite one-dimensional
$T_0$ space.

\begin{thm}\label{t:openorclosed}
Let $X$ be a finite $T_0$ space. Then $\dim X\leq 1$ if and only
if every singleton in $X$ is either open or closed.
\end{thm}
\begin{proof}
Let $X$ be a finite $T_0$ space. By
Proposition~\ref{p:tswisolated}, $X$ has an open point say $x_0$
and $U_{x_0}= \{x_0\}$. Since $\dim X\leq 1$, then $\dim \partial
(\{x_0\})= 0$ and so $\partial (\{x_0\})$ is discrete. This means
that each $y_0\in \partial (\{x_0\})$ is closed in $\partial
(\{x_0\})$. Since $\partial (\{x_0\})$ is closed in $X$, then
$\{y_0\}$ is also closed in $X$. Take $X^{'}= X\setminus
\mbox{cl}(\{x_0\})$ which is an open finite $T_0$ subspace of $X$.
By Proposition~\ref{p:openisolated}, $X^{'}$ has an open point set
$x_1$ and $U_{x_1}= \{x_1\}$ which is also open in $X$. Also $\dim
\partial (\{x_1\})= 0$, then $\partial (\{x_1\})$ is discrete. So,
each $y_1 \in \partial (\{x_1\})$ is closed in $\partial
(\{x_1\})$. Then $\{y_1\}$ is closed in $X$. Take $X^{''}=
X\setminus \mbox{cl}(\{x_0, x_1\})$. By continuing, the proof is
completed in one direction. Conversely, Suppose that singletons
are open or closed. By Lemma~\ref{l:clopendiscrete}, we have in
either cases discrete subspaces of $X$. Therefore the dimension of
each subspace is $X\leq 1$. Hence $\dim X\leq 1$.
\end{proof}

\begin{defn}
The height of a partially ordered set $(X, \leq)$ is the degree of
a longest element of the increasing sequence $x_1< x_2<....< x_n$
of elements of $X$.
\end{defn}

\begin{exam}
The discrete space has a height 1, since there are no comparable
elements i.e. no $x< y$. The space in which each point is either
open or closed of height 2, for $x< y$ is satisfied only for a
closed point $x$ and an open point $y$.
\end{exam}

Theorem~\ref{t:openorclosed} can be reformulate as follows: The
dimension of a finite space $X$ is the height of a partially
ordered set $(X, \leq)$ minus 1.

\begin{rem}
Recall that a space $(X, \tau)$ is called $T_{\frac{1}{2}}$ if
every generalized closed subset of $X$ is closed~\cite{levine2} or
equivalently if every singleton is either open or
closed~\cite{dunham}.\\
Observe that by Theorem~\ref{t:openorclosed}, every finite
$T_{\frac{1}{2}}$ space is a one dimensional $T_0$. It is also
clear that a space $X$ is $T_{\frac{1}{2}}$ if its non-closed
singletons are isolated.
\end{rem}

\section{Density in Finite Spaces}

In this section, we introduce some characterizations of density in
finite spaces.

\begin{prop}\label{p:openpdense}
In finite $T_0$ spaces, the set of all open points is dense.
\end{prop}
\begin{proof}
Let $X$ be a finite $T_0$. By Proposition~\ref{p:tswisolated}, $X$
contains an isolated point. If $O= \{x: \{x\}\mbox{ is open }\}$.
Then by Proposition~\ref{p:openisolated}, every open set $G$
containing $x$ must intersect $O$. This means that $x\in \overline
O$. Hence $O$ is dense.
\end{proof}

\begin{rem}
The condition $T_0$ is necessary in
Proposition~\ref{p:openpdense}. The singleton $\{x\}$ in $\tau_6$
of example~\ref{ex:thptnbd} is open but not dense.
\end{rem}

\begin{prop}\label{p:denseiffopen}
A subset $A$ is dense in a finite $T_0$ space $X$ if and only if
it contains all open points of $X$.
\end{prop}
\begin{proof}
Let $O=\{x:\{x\}=U_x\}$ be the set of open points in $X$. If $A$
is a dense subset of $X$, by density, then $A$ must intersect
$U_x$. So, it must contain $x$. The converse case follows readily
from~Proposition~\ref{p:openpdense}.
\end{proof}

\begin{cor}\label{c:denseconopen}
In finite $T_0$ spaces, each dense set contains an open dense set.
\end{cor}
\begin{proof}
Obvious from Proposition~\ref{p:openpdense} and
Proposition~\ref{p:denseiffopen}.
\end{proof}

In Corollary~\ref{c:denseconopen}, The condition $T_0$ is
necessary.

\begin{exam}
Consider the topology $\tau _6$ in Example~\ref{ex:thptnbd}. A
subset $\{a, c\}$ is dense. But it does not contain an open dense
subset.
\end{exam}

Some topological concepts induced by density. Recall that $A$ is
said to be codense (resp. nowhere dense, dense-in-itself) if
$\mbox { int }(A)= \phi$ (resp. $\mbox { int }\overline A= \phi$,
$A= \mbox{ d}(A)$), where $\mbox{ d}(A)$ denotes the set of
accumulation points of $A$. Dense-in-itself of any $A\subseteq X$,
equivalently that $A$ does not have any isolated points. We state
these notions in the finite spaces.

\begin{prop}
For arbitrary finite topological space $X$, the following are true
for any subset $A$ of $X$,
\begin{enumerate}
\item[(i)] It is co-dense if there is no upper bounds
belongs to $A$.
\item[(ii)] It is nowhere dense if there is no element has a
maximal element in $A$.
\item[(iii)] It is dense in-itself if $A$ contains all of its
lower bounds.
\end{enumerate}
\end{prop}
\begin{proof}
(i) and (ii) are obvious from Proposition~\ref{p:interiorset} and
Proposition~\ref{p:intcl} respectively.\newline (iii) Clear by
Proposition~\ref{p:closureset} and using the equality $A=
\overline A$.
\end{proof}

Recall that a topological space $(X, \tau)$ is said to be a
submaximal if each of its dense subsets is open.

\begin{prop}
Any finite $T_0$ space $X$ is submaximal if and only if it
contains at most a non-isolated point.
\end{prop}
\begin{proof}
The set of open points is dense by Proposition~\ref{p:openpdense}.
Also by the submaximality of $X$, every dense subset is open. Then
for every points $x, y, z \in X$ such that $x\leq y\leq z$,
$X\setminus \{y\}$ is not open. Therefore $y$ is not isolated.
Conversely, if all the points of $X$ are isolated, then the only
dense set is $X$ itself. This means that $X$ is submaximal.
Otherwise, if $X$ contains only non-isolated point $y$. Let $A$ be
any dense subset of $X$. By Proposition~\ref{p:denseiffopen}, $A$
contains all isolated points of $X$.
\end{proof}

\begin{cor}
Every finite $T_0$ submaximal spaces is $T_{\frac {1}{2}}$. It is
also one-dimensional space.
\end{cor}
\begin{proof}
It is a consequence Theorem~\ref{t:openorclosed}.
\end{proof}

\section{Some Weaker Forms of Continuity in Finite Spaces}

Continuous functions play an important role in topology. The
following proposition shows a corresponding definition of
continuity between finite spaces using the minimal neighborhood of
each point and its image.

\begin{prop}
Let $X$ and $Y$ be finite topological spaces, then the function
$f:X\longrightarrow Y$ is continuous at $x$ if and only if
$f(U_x)\subseteq U_{f(x)}$.
\end{prop}
\begin{proof}
Let $f:X\longrightarrow Y$ be a continuous function. Fix a point
$x\in X$. Since $U_{f(x)}$ is an open neighborhood of $f(x)$, then
by the continuity of $f$, $f^{-1} (U_{f(x)})$ is an open
neighborhood of $x$ and so $f(U_x) \subseteq U_{f(x)}$.
Conversely, let $W$ be an open set containing $f(x)$. By
assumption, $f(U _x)\subseteq U_{f(x)}\subseteq W$. Take $U_x= U$,
then $f(U)\subseteq W$ which shows that $f$ is continuous at $x$.
\end{proof}

Recall that a function $f: X\longrightarrow Y$ from a topological
space $X$ into a topological space $Y$ is precontinuous
\cite{mash} (resp. semicontinuous \cite{levine3},
$\gamma$-continuous \cite{abdel}) if the inverse image of each
open set in $Y$ is preopen (resp. semiopen, $\gamma$-open) in $X$.

In an arbitrary topological space, there is no connection between
precontinuity and semicontinuity. Each of them implies
$\gamma$-continuity. In the finite case, by
Remark~\ref{r:semipre}, every precontinuous function is
semicontinuous which coincides with $\gamma$-continuous. The
following implications show a new connection between continuity
and some kinds of near continuity.

$$\begin{array}{crccr}
\textbf{\large Continuity}&\Longrightarrow& \textbf {\large
Precontinuity}& \\
&&\Downarrow&& \\ 
\textbf{\large $\gamma $-continuity }&\Longleftrightarrow
&\textbf{\large Semicontinuity}&&
\end{array}$$

The following example shows that the converse is not always true.

\begin{exam}
Let $X= Y= \{x, y, z\}$ as in Example~\ref{ex:thptnbd}. The
mapping $f: (X, \tau_3)\longrightarrow (Y, \tau_5)$ which defined
by $f(x)= x$, $f(y)= z$ and $f(z)=y$ is $\gamma$-continuous, but
not precontinuous.
\end{exam}

A function $f:X\longrightarrow Y$ from a topological space $X$
into a topological space $Y$ is called preopen \cite{mash} (resp.
semiopen\cite{biswas}, $\gamma$-open\cite{abdel}) if the image of
each open set in $X$ is preopen (resp. semiopen, $\gamma$-open).

In arbitrary topological space, there is no connection between
preopen and semiopen functions. In a finite case, by
Remark~\ref{r:semipre}, a preopen function is semiopen which
coincides with $\gamma$-open. Each of them belongs to
$\gamma$-open function. The implications between these types of
functions and other corresponding ones are given by the following
diagram:

$$\begin{array}{crccr}
\textbf{\large Open function}&\Longrightarrow& \textbf {\large
Preopen function}& \\
&&\Downarrow&& \\ 
\textbf{\large $\gamma $-open function}&\Longleftrightarrow
&\textbf{\large Semiopen function}&&
\end{array}$$

The converse of these implications are not true, in general, as
the following example illustrates.

\begin{exam}
Let $X= \{x, y, z\}$ with a topology $\tau_5$ as in
Example~\ref{ex:thptnbd} and $Y= \{a, b, c, d\}$ with the topology
of minimal neighborhoods $U_a= \{a\}$, $U_b= \{b\}$, $U_c= U_d=
Y$. The function $f:X\longrightarrow Y$ is $\gamma$-open but not
preopen.
\end{exam}

\section{Continuity of Multifunctions in Finite Topological Spaces}

Let $X$ and $Y$ be two nonempty sets and $P(Y)$ be the collection
of all subsets of $Y$. Recall that the function $F:X
\longrightarrow P(Y)$ is called a multifunction~(\cite{cain},
p.186).

\begin{defn}
A function $F:X \longrightarrow P(Y)$ from a topological space $X$
into a topological space $Y$ is called:
\begin{enumerate}
\item \textit{upper semicontinuous mulitifunction} (abb.U.S.C) at a point
$x_0 \in X$ if for every open set $V$ in $Y$ such that $F(x_0)
\subset V$, there exists an open set $U$ containing $x_0$ such
that $F(U)\subset V$.
\item \textit{lower semicontinuous multifunction} (abb.L.S.C)at a point
$x_0 \in X$ if for every open set $V$ in $Y$ such that $F(x_0)
\cap V \neq \phi$, there exists an open set containing $x_0$ such
that $F(x) \cap V \neq \phi$ for each $x \in U$.
\end{enumerate}
\end{defn}

\begin{rem}
U.S.C is not necessarily L.S.C, and L.S.C is not necessarily
U.S.C.
\end{rem}

Now, the following theorems and examples give the relation beteen
the usual definition of continuity for single valued functions and
upper and lower semicontinuity for multifunctions.

\begin{defn}\cite{khal}
A space $X$ is said to be a connected ordered topological space
(abb. COTS) if for every three points subset $Y$ in $X$, there
exists $y\in Y$ such that $Y$ meets two connected components of
$X\setminus\{y\}$. In other words, for any three points one of
them separates the other two.
\end{defn}

\begin{rem}\label{r:single}
Let $X= [0,1]$ with usual topology and $f$ be a continuous
function from $X$ into itself. $\pi :X\longrightarrow Y$ is a
quotient function from $X$ into a finite COTS $Y$. The function
$g$ from $Y$ into itself is defined by $g(y)= \pi f
{\pi}^{-1}\{y\}$ is continuous. The following examples discuss the
continuity of $g$ in the case of single valued functions.
\end{rem}

\begin{exam}\label{e:pi}
Let $X= [0, 1]$ and $f:X \longrightarrow X$ be a continuous
function defined by \[f(x)= \left\{\begin{array}{r@{\quad:\quad}l}
x+ 1/2& 0 \leq x \leq 1/2\\1 & 1/2\leq x \leq 1 \end{array}
\right.\] and $Y= \{a, b, c, d, e\}$ is a COTS. $\pi$ is a
quotient function defined by \[\pi
(x)= \left\{\begin{array}{r@{\quad:\quad}l} a &x= 0\\b & 0<x<1/2\\
c& x= 1/2\\d &1/2<x<1\\ e &x= 1 \end{array} \right.\] Since $g(y)=
\pi f {\pi}^{-1}{(y)}$ for each $y\in Y$, so $g(a)= c$, $g(b)= d$,
$g(c)= e$, $g(d)= e$, $g(e)= e$. Then $g$ is continuous at each
point $y\in Y$.
\end{exam}

The following theorem is considered as an equivalent definition
for U.S.C and L.S.C in finite topological spaces with respect to
the minimal neighborhood of each point.

\begin{thm}\label{t:multif}
For finite topological spaces $X$ and $Y$ and multifunction
$F:X\longrightarrow P(Y)$, we have:
\begin{enumerate}
\item[(i)] $F$ is U.S.C at $x$ if and only if $F(\mathcal{U}_x)\subset 
\bigcup\limits_{y\in F(x)}\mathcal{U}_y$.
\item[(ii)] $F$ is L.S.C at $x$ if and only if $F(x')\cap \mathcal{U}_y\neq
\phi$ for all $x' \in \mathcal{U}_x$ and $y\in F(x)$.
\end{enumerate}
\end{thm}
\begin{proof}
(i). 
Let $F$ be U.S.C at $x\in X$. 
Then for every open set $V$such that $F(x)\subset V$. 
There exits an open set $U$ with $x\in U$ such that $F(U)\subset V$. 
Since $x\in U$, then $\mathcal{U} _x \subset \mathcal{U} $.
By assumption, we obtain $F(\mathcal{U} _x)\subset V$. 
Take 
$\textstyle{V= \bigcup\limits_{y\in F(x)}\mathcal{U} _y}$, 
then 
$F(\mathcal{U} _x)\subset \bigcup\limits_{y\in F(x)}\mathcal{U} _y$. 
Conversely, let $F(\mathcal{U} _x)\subset \bigcup\limits_{y\in F(x)}\mathcal{U} _y$. 
For $x\in \mathcal{U} _x$, take any $V\supset F(x)$. 
As $\mathcal{U} $, we choose $\mathcal{U} _x$, then 
$F(\mathcal{U} _x)\subset \bigcup\limits_{y\in F(x)}\mathcal{U} _y \subset V$. 
Therefore $F(\mathcal{U}_x) \subset V$ shows that $F$ is U.S.C. 

(ii). 
Let $F$ be L.S.C\ at $x\in X$. 
That means for every open set $V$ such that $F(x)\cap V\neq \phi$, there
is an open set $U$ with $x\in U$ such that $F(x')\cap V\neq \phi$ for all
$x'\in U$. 
Applying for a special $V$, take $V= \mathcal{U}_y$ and $y\in F(x)$. 
Since $\mathcal{U}_x\subset \mathcal{U}$, then $F(x')\cap \mathcal{U}_y \neq \phi$
for all $x'\in \mathcal{U}_x$ and $y\in F(x)$. 
Conversely, let $V$ be an arbitrary open set with 
$F(x)\cap V\neq \phi$. 
Take $y\in F(x)\cap V$, then $\mathcal{U}y \subset V$. 
As $\mathcal{U}$, we choose $\mathcal{U}_x$ and by assumption 
$\phi \neq F(x')\cap \mathcal{U}_y\subset F(x')\cap V$. 
Then $F(x')\cap V\neq \phi$ for all $x'\in \mathcal{U}_x$ and $F$ is L.S.C.
\end{proof}

U.S.C may not be L.S.C. 
It can be satisfied from the following example.

\begin{exam}
Let $X= Y= \{a, b\}$. Their minimal neighborhoods are $\mathcal{U}_a=
\{a\}$ and $\mathcal{U}_b= \{a, b\}$. The multifunction $F:
X\longrightarrow P(Y)$ is defined by $F(a)= \{a, b\}, F(b)=
\{a\}$. Then $F$ is L.S.C. It is not U.S.C since we have 
$F(\mathcal{U}_b)= \{a, b\}\not\subset \bigcup \limits_{y\in F(b)}\mathcal{U}_y=
\{a\}$ by Theorem~\ref{t:multif}(i), $\mathcal{U}_b=\{a, b\}, F(b)=
\{a\}$. Therefore $F$ is not U.S.C. If the multifunction
$G:X\longrightarrow P(Y)$ is defined by $G(a)= \{a\}, G(b)= \{a,
b\}$, then $G$ is U.S.C. It is not L.S.C at the point $b$. Since
by Theorem~\ref{t:multif}(ii), we have $\mathcal{U}_b= \{a, b\},
F(b)= \{a, b\}$ and $F(a)\cap \mathcal{U}_a= \{b\}\cap \{a\}= \phi$.
\end{exam}

The following examples discuss the upper and lower semi-continuity
of $g$.

\begin{exam}\label{e:pi1}
Let $X= [0, 1]$ with usual topology. The function $f$ from $X$
into itself defined by
\[f(x)= \left\{\begin{array}{r@{\quad:\quad}l}
3/4-2x& 0 \leq x < 1/4\\1/3x+1/6 & 1/4\leq x \leq 1 \end{array}
\right.\] and the quotient function $\pi: X\longrightarrow Y$,
where $Y$ is a finite COTS with five points defined by
\[\pi (x)= \left\{\begin{array}{r@{\quad:\quad}l} a &x= 0\\
b &x= 0<x<1/2\\c& x= 1/2\\d &1/2<x<1\\ e &x= 1
\end{array} \right.\]
The multifunction $g: Y\longrightarrow P(Y)$ is defined by $g(y)=
\pi f {\pi}^{-1}{y}$. Then $ g(a)= d, g(b)=\{b, c,d \}, g(c)= b,
g(d)= b, g(e)=c$. $g$ is L.S.C, but not U.S.C since the point $y=
a$, $\mathcal{U}_a= \{a, b\}$, $g(a)= \{d\}$ and $g(\mathcal{U}_a)= \{b,
c, d\} \not\subset \bigcup\limits_{y\in g(a)}\mathcal{U}_y)= \{d\}$.
\end{exam}

Now, we consider some properties of U.S.C and L.S.C for
multifunctions. These make a connection between multifunctions and
dimension theory.

\begin{thm}\label{t:contdim}
Let $\pi:X\longrightarrow Y$ be a continuous function from a
topological space $X$ into a one-dimensional $T_0$ space $Y$. For
any closed point in $Y$, ${\pi}^{-1}(y)$ is only one point and for
any open point $z\in \mathcal{U}_y$, ${\pi}^{-1}(y)\in
cl({\pi}^{-1}(z))$, then $\pi$ is an open function.
\end{thm}
\begin{proof}
Let $W$ be an open set in $X$. Our aim is to prove that $\pi (W)=
\bigcup\limits_{y\in \pi (W)}\mathcal{U}_y$. It is clear that $\pi
(W)\subset \bigcup\limits_{y\in \pi (W)}\mathcal{U}_y$. So, it is
enough to prove that $\mathcal{U}_y \subset \pi (W)$ for each $y\in
\pi (W)$. There are two cases: If $y$ is an open point, then
$\mathcal{U}_y= \{y\}\subset \pi (W)$ and we are done. If $y$ is a
closed point, by the assumption, ${\pi}^{-1}(y)= x$. Let $z$ be an
open point such that $z\in \mathcal{U}_(y)$, by hypothesis, $x\in
cl({\pi}^{-1}(z))$. Since $W$ is an open neighborhood of $x$, then
${\pi}^{-1}(z)\cap W \neq \phi$. For a point $v\in
{\pi}^{-1}(z)\cap W \neq \phi$, we have $z= \pi(v)\in \pi(W)$.
Therefore, $\mathcal{U}_y\subset \pi (W)$. This completes the proof.
\end{proof}

Under the previous assumptions of Theorem~\ref{t:contdim}, $\pi$
is not a closed function.

\begin{exam}
Let $X=[0, 1]$ with usual topology and $Y$ is a COTS with three
points. The function $\pi:X\longrightarrow Y$ is
defined by \[\pi (x)= \left\{\begin{array}{r@{\quad:\quad}l} a &x= 0\\
b &0<x<1\\c& x= 1 \end{array} \right.\] The assumption of Theorem
\ref{t:contdim} hold. $\pi$ is not closed because if we take a
closed set $F\subset (0, 1)$, then $\pi (F)=\{b\}$ is not closed.
\end{exam}

\begin{prop}
For a topological space $X$ and one-dimensional $T_0$-space $Y$,
consider a continuous function $f$ and a quotient a function $\pi$
from $X$ into $Y$. If the multifunction $g:Y\longrightarrow P(Y)$
is defined by $g= \pi f {\pi}^{-1}$, then $g(y)$ is a singleton
for each closed point $y\in Y$.
\end{prop}
\begin{proof}
Obvious, since $f$ and $\pi$ are single -valued functions.
\end{proof}

\begin{thm}\label{t:condit}
For a space $X$ and one-dimensional $T_0$-space $Y$, assume the
following assumptions hold.
\begin{enumerate}
\item[(1)] The function $\pi:X\longrightarrow Y$ is continuous.
\item[(2)] For every closed point $y$, ${\pi}^{-1}(y)$ is one
point.
\item[(3)] For every open point $z\in \mathcal{U}_y$, ${\pi}^{-1}(z)=
x$.
\item[(4)] A multifunction $g:Y\longrightarrow P(Y)$ defined by
$g= \pi f{\pi}^{-1}$.
\end{enumerate}
Then $g$ is L.S.C.
\end{thm}
\begin{proof}
Let $y_0\in Y$ and $V$ be an open set in $Y$ such that $g(y_0)\cap
V\neq \phi$. If $y_0$ is an open point, then we can take 
$\mathcal{U}_{y_0}=\{y_0\}$ and so $ g$ is L.S.C. If $y_0$ is a closed
point, by assumption, for any open point $y\in \mathcal{U}_{y_0}$,
${\pi}^{-1}(y_0)= x\in cl({\pi}^{-1}(y))$. By Proposition 30,
$g(y_0)$ is a singleton and so $f^{-1}{\pi}^{-1}(V)$ is an open
neighborhood of $x$. Therefore ${\pi}{-1}(V)\cap W \neq \phi$ for
all $y\in \mathcal{U}_{y_0}$ and $W= f^{-1}{\pi}^{-1}(V)$. This shows
that $g$ is L.S.C.
\end{proof}

Under the conditions of Theorem~\ref{t:condit}, the multifunction
is not necessarily to be U.S.C.

\begin{exam}
Let $X$ and $Y$ and the quotient function as in
Example~\ref{e:pi1}. The function $f$ is defined by \[f(x)=
\left\{\begin{array}{r@{\quad:\quad}l} 2x-1& 0 < x < 1/2\\1-2x &
1/2< x < 1\\1&x\in \{0, 1\} \end{array} \right.\] We have $ g(a)=
e, g(b)=\{b, c,d \}, g(c)= a, g(d)= \{b, c,d\}, g(e)=e$. $g$ is
not U.S.C since at the point $y_0= a$, $\mathcal{U}_a= \{a, b\}$,
$g(a)= \{e\}$, we have $g(\mathcal{U}_a)= \{b, c, d, e\} \not\subset
\bigcup\limits_{y\in g(a)}\mathcal{U}_y)= \{d, e\}$.
\end{exam}

\subsection*{Acknowledgment}
The results of this paper were obtained during my Ph.D. studies at Tanta 
University and are also contained in my thesis.\footnote{A short version 
of the thesis, supervised by C. Bandt, was distributed as a preprint at 
Ernst-Moritz-Arndt-Universit\"at Greifswald \cite{elatik}.} 
I would like to express my deep gratitude to my supervisors whose guidance
and support were crucial for the successful completion of this paper.

\providecommand{\bysame}{\leavevmode\hbox to3em{\hrulefill}\thinspace}
\providecommand{\MR}{\relax\ifhmode\unskip\space\fi MR }
\providecommand{\MRhref}[2]{%
  \href{http://www.ams.org/mathscinet-getitem?mr=#1}{#2}
}
\providecommand{\href}[2]{#2}

\end{document}